\documentclass{article}
\usepackage{graphicx} 
\usepackage{amsmath}
\usepackage{amsthm}
\newtheorem{theorem}{Theorem}[section]
\newtheorem{definition}{Definition}[section]
\newtheorem{lemma}[theorem]{Lemma}
\usepackage[dvipsnames]{xcolor}
\usepackage{imakeidx}
\makeindex
\usepackage{amssymb}
\usepackage{tikz}
\usetikzlibrary{arrows.meta}
\usepackage{hyperref}
\usepackage{natbib}
\title{Hausdorff Measure and Dimension with Examples}
\author{Umberto Michelucci\\{ \url{umberto.michelucci@hslu.ch}}}
\date{November 2025}

\begin{document}

\maketitle

\begin{abstract}
This document offers a concise introduction to the mathematical theory and practical application of the Hausdorff Measure and Dimension. The primary objective is to clarify and rigorously detail the two most common methods used for calculating the dimension of a set, ensuring all calculation details are transparent for the reader.

The paper first establishes the theoretical groundwork by reviewing the definitions of the Hausdorff measure, proving the dimensional invariance under changes to the shape of the covering sets, and confirming the dimensional property of open sets. It then introduces the two main methodologies.

The first is the Lower and Upper Bound Estimation, which uses the  relationship between the measure $\mathcal{H}^s(A)$ and the dimension $\dim_{\mathcal{H}}(A)$. This method emphasizes the use of the Mass Distribution Principle for establishing the lower bound, which is essential when the Lebesgue measure of the set is zero. The second, more computationally efficient method is the Similarity Dimension Method, introduced via the definitions of Similitudes, Iterated Function Systems (IFS), and the Moran-Hutchinson Theorem.

Both methodologies are applied rigorously to classic examples: the Unit Square (yielding dimension $2$) and the Cantor Set (yielding the fractional dimension $\log 2 / \log 3$). The paper serves as a detailed, step-by-step guide intended to make the application of these fundamental fractal geometry concepts clearer and easier to understand.
\end{abstract}

\newpage
\tableofcontents
\newpage

\section{Introduction}

These short notes, covers an introduction on two methods on how to calculate the Hausdorff dimension: the lower and upper bound approach (using the mass distribution principle) and the similarity method. I apply the methods to two sets: a unit square and the Cantor set. I hope to make the application of the methods clearer for the reader by giving all details of the calculations in the most rigorous way possible.

\section{Definitions and Properties}
The Hausdorff measure\index{Hausdorff measure} can be defined in several way. For example in \citep{lio_analysis_nodate} it is defined by the following.
\begin{definition}[\textcolor{red}{Spherical Hausdorff Measure}]
\index{Hausdorff measure}
    For $s\geq 0$, $\delta>0$ and $\emptyset \neq A \subseteq \mathbb{R}^n$ we define
    \begin{equation}
        \mathcal{H}_\delta^s(A) = \inf \left\{
            \sum_{k\in I} r_k^s , A\subseteq\bigcup_{k\in I} B(x_k, r_k), 0<r_k<\delta
        \right\}
    \end{equation}
    we also set $\mathcal{H}_\delta^0(\emptyset)=0.$ The set of indices $I$ is at most countable (so $I\subseteq \mathbb{N}$).
\end{definition}
Note how this definition uses balls\index{ball} as a covering\index{covering}. Furthermore, it is defined such that $r_k>0$, so that $r_k$ cannot be zero. in this case, the set of all balls is \textbf{not} technically a covering (according to the definition \ref{eq:covering} from \citep{lio_analysis_nodate}), since, for example, $\emptyset$ is not included in the set of balls (see definition \ref{eq:covering}). Some authors use the less restrictive definition of a covering (given in \ref{def:covering2}) so in this sense, the set of balls is a covering. A good question is why we use balls in the definitions and not generic sets. Indeed, this is what other authors use. For example, in \citep{falconer2003fractal} the Hausdorff dimension is defined as the following
\begin{definition}[\textcolor{red}{Hausdorff Definition (Generic Sets)}]
\label{def:haus-generic-sets}
    Let $F\subseteq \mathbb{R}^n$, $s>0$ with $s\in \mathbb{R}$. For any $\delta>0$ we define
    \begin{equation}
        \mathcal{H}_\delta^s(A) = \inf \left\{
            \sum_{i=1}^\infty |U_i|^s \ | \ {U_i} \text{ is a $\delta$-cover of F}
        \right\}
    \end{equation}
    Thus we are interested at all covers of $F$ by sets of diameters at most $\delta$. Where $|U_i|$ indicates the diameter of the set $U_i$ (as defined in \citep{falconer2003fractal}).
\end{definition}
Note that while in the definitions with balls the radius must be at most $\delta$, in the generic one it is the diameter of the sets that must be at most $\delta$. Note how there is a factor of 2 at play. Note also that naturally, those different definitions give different \textit{measure} values of sets, but the useful property is that the \textbf{Hausdorff dimension} (that we define in Definition \ref{def:hausdorff-dim}) will be the same regardless of the shape of sets that we decide to use in the definition as long as we use the diameter in the definition.

I want to cite two important theorems of which I will give no proof (proofs can be found in \citep{lio_analysis_nodate}).

\begin{theorem}[\textcolor{red}{$\mathcal{H}_\delta^s$ is a measure on $\mathbb{R}^n$}]\index{Hausdorff measure!measure}
$\mathcal{H}_\delta^s$ is a measure on $\mathbb{R}^n$.
\end{theorem}

\begin{theorem}[\textcolor{red}{$\mathcal{H}_\delta^s$ is Borel Regular Measure on $\mathbb{R}^n$}]\index{Hausdorff measure!Borel regular}
$\mathcal{H}_\delta^s$ is Borel Regular Measure on $\mathbb{R}^n$
\end{theorem}


\subsection{Hausdorff Dimension}
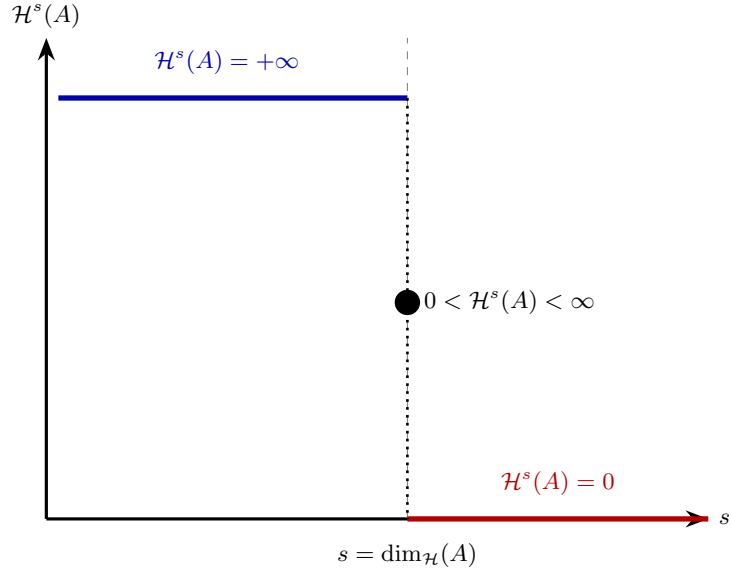
\begin{figure}[h]
    \centering
    \begin{tikzpicture}[
        >={Stealth[length=3mm]},
        scale=1.6, 
        font=\small
    ]

    \def\dimH{3.0}
    \def\maxy{4.0}
    
    \draw[->, very thick] (0,0) -- (5.5,0) node[right] {$s$};
    \draw[->, very thick] (0,0) -- (0,\maxy) node[above] {$\mathcal{H}^s(A)$};

    \draw[dashed, gray] (\dimH, 0) -- (\dimH, \maxy);
    \node at (\dimH, -0.3) {$s = \dim_{\mathcal{H}}(A)$};
    
    
    \draw[line width=2pt, blue!70!black] (0.1, 3.5) -- (\dimH, 3.5);
    \node[blue!70!black] at (1.5, 3.8) {$\mathcal{H}^s(A) = +\infty$};
    
    \draw[line width=2pt, red!70!black] (\dimH, 0) -- (5.5, 0);
    \node[red!70!black] at (4.25, 0.3) {$\mathcal{H}^s(A) = 0$};
    
    
    \draw[line width=1pt, dotted, black] (\dimH, 3.5) -- (\dimH, 0);
    
    \fill[black] (\dimH, 1.8) circle (3pt);
    \node[right, shift={(0.1, 0)}] at (\dimH, 1.8) {$0 < \mathcal{H}^{s}(A) < \infty$};

    \node[anchor=north] at (1.5, -0.3) {};
    \node[anchor=north] at (4.25, -0.3) {};

    \end{tikzpicture}
    \caption{The Phase Transition of the Hausdorff Measure $\mathcal{H}^s(A)$. The Hausdorff Dimension is the critical exponent ($s$) where the measure drops from infinity to zero.}
    \label{fig:hausdorff_switch_clean}
\end{figure}

To define the Hausdorff dimension, we first need another definition. Let us define $\mathcal{H}^s$ as
\begin{definition}[\textcolor{red}{$s$-dimensional Hausdorff measure}]
\begin{equation}
    \mathcal{H}^s = \lim_{\delta \downarrow 0} \mathcal{H}^s_\delta
\end{equation}
\end{definition}
\noindent Note that in general $\mathcal{H}^s$ can also be zero or infinite. Here are some properties of $\mathcal{H}^s$.
\begin{lemma}[\textcolor{red}{Properties of $\mathcal{H}^s$}]
\label{eq:lemma1}
Here are some useful properties that we report with proofs (for proofs see \citep{lio_analysis_nodate}).
    \begin{enumerate}
        \item Let $A\subseteq\mathbb{R}^n$ and $0\leq s<t<+\infty$. Then
        \begin{enumerate}
            \item $\mathcal{H}^s(A) < +\infty \Rightarrow \mathcal{H}^t(A)=0$
            \item $\mathcal{H}^t(A) >0 \Rightarrow \mathcal{H}^s(A)=+\infty$
        \end{enumerate}
        \item for $A\subseteq \mathbb{R}^n$, if $s>n$ then $\mathcal{H}^s(A)=0$.
    \end{enumerate}
\end{lemma}
\begin{definition}[\textcolor{red}{Hausdorff Dimension}]
\label{def:hausdorff-dim}
The Hausdorff dimension\index{Hausdorff dimension} of a set $A$ can be defined by\index{Hausdorff dimension}
\begin{equation}
    \dim_\mathcal{H}(A)= \inf \{s\geq 0 : \mathcal{H}^s(A)=0 \} = \sup \{s\geq 0 : \mathcal{H}^s(A)=+\infty \}
\end{equation}
\end{definition}
\noindent The behaviour of the Hausdorff dimension is visually depicted in Figure \ref{fig:hausdorff_switch_clean}.
\noindent An important property of the Hausdorff dimension is the following.
\begin{lemma}
    If $A\subseteq \mathbb{R}^n$ is a non-empty open set, then $\dim_\mathcal{H}(A)=n$.
\end{lemma}
Another important property, and one that will help us to calculate the Hausdorff dimension in an efficient way, is the following.
\begin{lemma}
The Hausdorff dimension is invariant\index{Hausdorff dimension!invariance} under scaling and changing the shape of the sets used in the definition. In particular, it is true that
\begin{equation}
    \dim_\mathcal{H}(\lambda A) = \dim_\mathcal{H}(A) \ \text{ for } \lambda>0, \lambda \in \mathbb{R}
\end{equation}
\end{lemma}
\noindent The proof is relatively easy but I will not report it here.

\subsection{Hausdorff Dimension Invariance}
As an example we can show how the Hausdorff dimension is not dependent on the shape of the sets used in the definition of $\mathcal{H}_\delta^s$. 
Consider a simple example for illustrative purposes.
 Recall definition \ref{def:haus-generic-sets}. We will consider the general definition with arbitrary sets ($\mathcal{H}_\delta^s$) and one where we consider only balls for the sets $U_i$ that we will indicate with $\mathcal{H}_\delta^{s,\circ}$. 
\begin{lemma}[\textcolor{red}{Invariance of Hausdorff Dimension from the Shape of Sets}]
Given our definition it holds that
\begin{equation}
    \operatorname{dim}_{\mathcal{H}}(A) = \operatorname{dim}_{\mathcal{H^{\circ}}}(A) 
\end{equation}
\end{lemma}
\begin{proof}
 
Now since clearly 
$\{\text{set of all possible balls}\}\subseteq \{ \text{all possible sets $U_i$} \}$ it follows that for an arbitrary set $E$
\begin{equation}
    \mathcal{H}^s(E)\leq \mathcal{H}^{s,\circ}(E)
\end{equation} since we are searching the infimum.
Now consider $x_i\in U_i$. $U_i$ is contained in the ball $B(x_i,r_i)$ with $r_i=\operatorname{diam} U_i$. Surely, $E\subseteq \bigcup_{i=1}^\infty B_i$. But each ball has $\operatorname{diam}B_i=2 \operatorname{diam}U_i$ so clearly
\begin{equation}
    \sum_i (\operatorname{diam} B_i)^s = 2^s  \sum_i (\operatorname{diam}U_i)^s
\end{equation}
thus
\begin{equation}
    \mathcal{H}^{s,\circ}\leq 2^s \mathcal{H}^s
\end{equation}
so we can say
\begin{equation}
\label{eq:dis1}
    \mathcal{H}^s(E) \leq \mathcal{H}^{s,\circ} (E)\leq 2^s \mathcal{H}^s(E)
\end{equation}
We can use Lemma \ref{eq:lemma1} and notice that $\mathcal{H}^s$ can assume only three values: 0, a finite value or $+\infty$. 
let us call $s=\dim_\mathcal{H}(E)$ and $s^\circ=\dim_\mathcal{H^{\circ}}(E)$.
\begin{enumerate}
    \item[(i)] If $t>s$ then $\mathcal{H}^{t}(E)=0$, but from Equation (\ref{eq:dis1}) it follows that $\mathcal{H}^{t,\circ}(A)\leq 2^t 0=0$. Thus $s^\circ \leq s$. $s$ is an upper bound for $s^\circ$.
    \item[(ii)] If $t<s$  then $\mathcal{H}^t(E)=+\infty$ and hence $\mathcal{H}^{s,\circ}(E) \geq +\infty$ and thus $\mathcal{H}^{s,\circ}(E)=+\infty$. So we conclude that $s^\circ \geq s$. So $s$ is a lower bound for $s^\circ$.
\end{enumerate}   
So from the two previous points we conclude $\dim_\mathcal{H^\circ}(E)=\dim_\mathcal{H}(E)$.
\end{proof}
\paragraph{Note} we will not discuss this further here, but for the interested reader, the Hausdorff dimension is also invariant under isometries (generally speaking transformation that preserve distance in a metric space). You can find a lot of interesting information in, for example, \citep{hao_hausdorff_nodate, cao_several_2023}.

\section{Methods for Calculating the Hausdorff Dimension}

There are many methods to calculate the Hausdorff dimension. A brief list, which is surely not complete, includes: (i) the upper and lower bound estimation method (discussed here in Section \ref{sec:lower-upper}), (ii) potential theory method (found in \citep{cao_several_2023}), (iii) project theorem for dimension estimation (found in \citep{cao_several_2023}), (iv) fractal multiplication for dimension estimation (found in \citep{cao_several_2023}), (v) similarity method (discussed here in Section \ref{sec:sim-meth}).

\subsection{Lower and Upper Bound (Mass Distribution Principle)}
\index{mass distribution principle}\index{Hausdorff dimension!lower bound}\index{Hausdorff dimension!upper bound}
\label{sec:lower-upper}

This method is based on the specific property of the Hausdorff Dimension (discussed in Lemma \ref{eq:lemma1}). In particular, if $s=\dim_\mathcal{H}(A)$, then $\mathcal{H}^t(A)=0$ for all $t>s$ and $\mathcal{H}^t(A)=+\infty$ for $t<s$. Keeping this in mind, to calculate the Hausdorff dimension we can use a two step approach.

\paragraph{Step 1 - Upper Bound}
\index{Hausdorff dimension!upper bound}

\noindent This step is relatively easy. If we find a specific covering such that $\sum_i |U_i|^s=0$ then we can say that $\dim_\mathcal{H}(A)\leq s$. Now, in general, this is not that useful in practice, as is stated as it might be quite difficult to achieve, since a finite covering with sets with finite diameter will always have $\sum_i |U_i|^s>0$. A good approach is to explicitly construct a sequence of $\delta$-coverings $\{U_{i}^{(\delta)}\}$ such that for a certain value of $s$ then $\lim_{\delta \to 0} \sum_i |U_i^{(\delta)}|^s = 0$. Now that means that for a given $s$ value $\mathcal{H}^s(A)=0$ and thus, the $\dim_\mathcal{H}(A)\leq s$ (due to how the Hausdorff dimension is defined as an infimum).

\paragraph{Step 2 - Lower Bound}
\index{Hausdorff dimension!lower bound}

\noindent This step is usually a bit more convoluted and involves a more complex process. The first step is about finding an efficient covering; the lower bound (this step) is essentially about proving that you cannot cover the set efficiently, even if you try. Since the definition of $\mathcal{H}^s(A)$ involves an $\inf$ over all possible coverings, you cannot just pick one bad covering and be done. You need a method that works against any covering.

The idea is based on the so-called \textbf{mass distribution principle}. This can be formulated in this way.
    Let $A$ be a non-empty Borel subset of $\mathbb{R}^n$ and $s \ge 0$.
If there exists a mass distribution $\mu$ on $A$ (i.e., a finite non-zero Borel measure such that $0 < \mu(A) < \infty$ and $\text{supp}(\mu) \subset A$) and constants $c > 0$ and $\delta_0 > 0$, such that the following scaling condition holds for every Borel set $U \subset \mathbb{R}^n$ with diameter $|U| < \delta_0$:
$$
\mu(U) \leq c \cdot |U|^s
$$
Then the $s$-dimensional Hausdorff measure of $A$ is bounded below by a positive constant:
$$
\mathcal{H}^s(A) \geq \frac{\mu(A)}{c} > 0
$$
Consequently, the Hausdorff dimension of $A$ is at least $s$:
$$
\dim_\mathcal{H}(A) \geq s
$$
The trick is naturally to find a mass distribution for each specific case. If interested, the proof of this can be found in \citep{hao_hausdorff_nodate}.

\subsection{Hausdorff Dimension of a Square}
\index{Hausdorff dimension!square}
Let us apply as a demonstration the lower and upper bound method described in the previous section to a square $S=[0,1]\times [0,1]\in \mathbb{R}^2$. Of course we would expect it to have a dimension equal to two.

\paragraph{Step 1 - Upper Bound.}
As we have discussed we need to find a sequence of $\delta$-covering such that 
\begin{equation}
    \lim_{\delta \to 0} \sum_i |U_i^{(\delta)}|^s = 0
\end{equation}
Such a sequence can be constructed by considering a set of smaller squares $U_i$ of side $\delta = L/N$. Their diameter is their diagonal, of length $\delta\sqrt{2}=(L/N)\sqrt{2}$ (the square diagonal). 
Now we can write (note that $\delta\rightarrow 0$ corresponds to $N\rightarrow +\infty$) the following.
\begin{equation} 
    \lim_{N \rightarrow +\infty}\sum_{\text{all squares}} (\text{diam } U_i)^s 
    = \lim_{N \rightarrow +\infty} N^2 \left( \sqrt{2}\frac{L}{N}\right)^s= \lim_{N \rightarrow +\infty}(\sqrt{2}L)^s N^{2-s}
\end{equation}
this is 0 for $s>2$. Thus, we can immediately say that $\dim_\mathcal{H}(S)\leq 2$, or in other words we have found an upper bound for the Hausdorff dimension.

\paragraph{Step 2 - Lower Bound.}

We apply the \textbf{Mass Distribution Principle} with the test dimension $s=2$.
We can define the measure $\mu$ on $A$ as the $2$-dimensional Lebesgue measure (area) restricted to the square S.
$$
\mu(E) = \lambda^2(E \cap S) \quad \text{for any Borel set } E \subset \mathbb{R}^2.
$$
We need to verify that our $\mu$ is a valid mass distribution. We need to check two things.
\begin{itemize}
    \item \textbf{Finite and Non-Zero Mass:} $\mu(S) = \lambda^2(S) = 1$. Thus, $0 < \mu(S) < \infty$.
    \item \textbf{Support on $S$:} By construction, $\mu(\mathbb{R}^2 \setminus S) = 0$.
\end{itemize}
So we conclude that $\mu$ is a valid mass distribution.
As a second step we need to verify the Scaling Condition $\mu(U) \leq c \cdot |U|^s$.
We must show that for any set $U$ with small diameter $|U|$, the mass $\mu(U)$ is bounded by $|U|^2$ times a constant $c$.

Let $U \subset \mathbb{R}^2$ be any Borel set. Since $U$ is contained in a closed ball $B$ of radius $r = |U|/2$, and the area of this ball is $\pi r^2$, we have:
$$
\mu(U) = \lambda^2(U \cap S) \leq \lambda^2(U) \leq \lambda^2(B)
$$
Substituting $r = |U|/2$ into the area formula for the ball:
$$
\mu(U) \leq \pi r^2 = \pi \left(\frac{|U|}{2}\right)^2 = \left(\frac{\pi}{4}\right) \cdot |U|^2
$$
We have successfully established the scaling condition:
$$
\mu(U) \leq c \cdot |U|^s \quad \text{with } s=2 \text{ and } c=\frac{\pi}{4}.
$$
Since the conditions of the principle are satisfied with $s=2$, we conclude:
$$
\mathcal{H}^2(S) \geq \frac{\mu(S)}{c} = \frac{1}{\pi/4} = \frac{4}{\pi} > 0
$$
\noindent Therefore, the Hausdorff dimension is bounded below by 2.
$$
\dim_\mathcal{H}(S) \geq 2
$$
Combining this with the upper bound $\dim_\mathcal{H}(S) \leq 2$ gives the expected result: $\dim_\mathcal{H}(S) = 2$.

\subsection{Hasudorff dimension of the Cantor Set}
\index{Hausdorff dimension!cantor set}\index{Cantor set!definition}

The Cantor Set, denoted $C$ is famous for being compact, totally disconnected, and having a non-integer Hausdorff dimension.
\begin{definition}[\textcolor{red}{Cantor Set}]

\noindent The set $C$ is formally defined as the intersection of a sequence of sets $\{C_n\}_{n=0}^{\infty}$:
$$
C = \bigcap_{n=0}^{\infty} C_n
$$
where the sequence $\{C_n\}$ is constructed iteratively in the following steps:

\begin{enumerate}
    \item \textbf{Stage 0 ($\mathbf{C_0}$):} Start with the unit interval:
    $$
    C_0 = [0, 1]
    $$

    \item \textbf{Stage 1 ($\mathbf{C_1}$):} Remove the open middle third $(\frac{1}{3}, \frac{2}{3})$ from $C_0$. This leaves $N=2$ closed intervals, each of length $r=1/3$:
    $$
    C_1 = \left[0, \frac{1}{3}\right] \cup \left[\frac{2}{3}, 1\right]
    $$

    \item \textbf{Stage $\mathbf{n}$ ($\mathbf{C_n}$):} If $C_{n-1}$ is the union of $2^{n-1}$ disjoint closed intervals, $C_n$ is formed by removing the open middle third from each of those $2^{n-1}$ intervals. $C_n$ thus consists of $2^n$ closed intervals, each of length $3^{-n}$:
    $$
    C_n = \bigcup_{k=1}^{2^{n}} I_{n,k}, \quad \text{where } |I_{n,k}| = \left(\frac{1}{3}\right)^n
    $$

\end{enumerate}
You can see a diagram explaining the Cantor set in Figure \ref{fig:cantor}.
\begin{figure}
    \centering
    \includegraphics[width=1\linewidth]{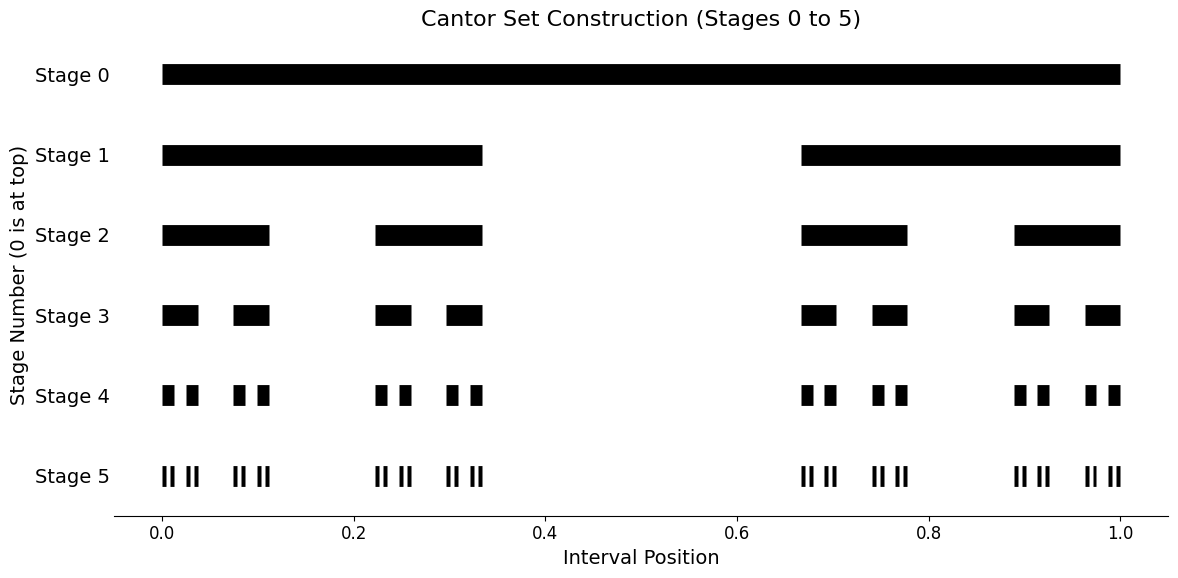}
\caption{The construction of the Middle-Thirds Cantor Set (C) over six stages (n=0 to n=5). At each stage, the open middle third is removed from every remaining interval. This process highlights the self-similar nature of the set, where the initial interval is scaled down by a factor r=1/3 and copied N=2 times at each iteration. This scaling relationship (N=2,r=1/3) is key to calculating the Hausdorff dimension using the Similarity Dimension method, yielding $\dim_\mathcal{H}(C)=\log 2/\log 3$ (see Section \ref{sec:sim-meth}).}
    \label{fig:cantor}
\end{figure}
    
\end{definition}
\noindent $C$ is therefore the set of all points in $[0, 1]$ that are \textbf{never removed} during this infinite process.
Let us calculate now its dimension using the lower/upper bound approach. Note that it is much easier to calculate its dimension with the similarity method described in Section \ref{sec:sim-meth} (go and check it out). The method with the lower/upper bounds requires more work.

\paragraph{Step 1 - Upper Bound.} As we have discussed and done for the square, we need to finda a sequence of $\delta$-covering such that
\begin{equation}
    \lim_{\delta \downarrow 0} \sum_i |U_i^{(\delta)}|^s = 0
\end{equation}
At each step $n$ (where we have the set $C_n$) We can consider balls in one dimension with diameter $(1/3)^n$. At iteration $n$ we have $2^n$ balls, therefore
\begin{equation}
    \sum_i (\operatorname{diam} U_i)^s = 2^n \left( \frac{1}{3}\right)^{ns}
\end{equation}
for $n\rightarrow \infty$ (since we want the dimension of the set $C$) this is going to zero for
\begin{equation}
    \frac{2}{3^s}<1
\end{equation}
and thus for 
\begin{equation}
    s>\frac{\log 2}{\log 3}
\end{equation}
so we have an upper bound for the Hausdorff dimension
\begin{equation}
    \dim_\mathcal{H}(C) \leq s = \frac{\log 2}{\log 3}.
\end{equation}

\paragraph{Step 2 - Lower Bound.} As we have done for the square let's apply the mass distribution principle with the test dimension $s = \log 2/\log 3$. We need to create or find an appropiate mass to use for the method. Clearly the Lebesgue measure $\mathcal{L}^1$ is not appropriate as we know that $\mathcal{L}^1(C)=0$. We can build an appropriate measure in this way. 
\begin{enumerate}
    \item \textbf{Step 0}: $\mu(C_0)=1$
    \item \textbf{Step 1}: the mass is divided equaly on the two intervals (indicated with $I_{1,1}, I_{1,2}$)so $\mu(I_{1,1})=\mu(I_{1,2})=1/2$.
    \item \textbf{Step $n$}: each interval $I_{n,k}$ gets $(1/2)^n$. 
\end{enumerate}
This will give us
\begin{equation}
    \mu(C)=\mu(\bigcap_{i=1}^\infty C_n) = \sum_{i=1}^{2^n}   \mu(I_{n,k})=2^n \frac{1}{2^n} = 1
\end{equation}
So we have for $n\rightarrow \infty$ $\mu(C)=1$.
Now we need to check that $\mu(U)\leq c |U|^s$ for our candidate $s$. Note as before we define for an arbitrary set $\mu(A)\equiv \mu (A\cap S)$. Our $U$ is an arbitrary set with diameter at most $|U|=\delta$.

Now comes the harder part. Suppose we are considering a generic set $U$ with diameter $\delta$. We can surely find an $n$ s.t.
\begin{equation}
    \left( \frac{1}{3} \right)^{n+1}\leq |U| \leq \left( \frac{1}{3} \right)^{n}
\end{equation}
Now at stage $n$ such a set, given its diameter, can only cover 
Such an interval can touch up to two intervals at their end point (remember the distance between intervals $I_{n,k}$ in the Cantor set at stage $n$ is $3^{-n}$.

Now recall also that each interval $I_{n,k}$ has a mass of $(1/2)^n$. 
If $U$ covers two intervals we will have in the worst case scenario
\begin{equation}
    \mu(U) = \mu(U\cap C_n) \leq 2\cdot 2^{-n}
\end{equation}
At the same time we have, by construction, $|U|\geq 3^{-(n+1)}$. Recall now that we are testing $s=\log 2/\log 3$ and thus $2=3^s$. So it holds
\begin{equation}
    \mu(U) \leq 2\cdot 2^{-n} = 2\cdot (3^s)^{-n} =2\cdot 3^s |U|^s = 4  |U|^s 
\end{equation}
this is valid for all $\delta$ (and consequently for all $n$).
Thus we can conclude that, in the limit for $n\rightarrow \infty$ we have
\begin{equation}
    \dim_\mathcal{H}(C) \geq s = \frac{\log 2}{\log 3}
\end{equation}
From the two steps we can then finally conclude that
\begin{equation}
    \dim_\mathcal{H}(C) = \frac{\log 2}{\log 3}
\end{equation}

\subsection{The Similarity Dimension Method}
\label{sec:sim-meth}
\index{similarity dimension method}

This method determines the Hausdorff dimension ($\dim_{\mathcal{H}}(A)$) by solving a characteristic equation based on the set's scaling properties, provided specific conditions are met.
To give a proper formal treatment of this method, some definitions are needed.


\begin{definition}[\textcolor{red}{Similitude (Similarity Map)}]
    \index{simmilarity map}
A function $f: \mathbb{R}^n \to \mathbb{R}^n$ is a \textbf{similitude} if there exists a constant $r \in (0, 1)$ such that for all $x, y \in \mathbb{R}^n$:
$$
|f(x) - f(y)| = r |x - y|
$$
The constant $r$ is the \textbf{scaling ratio} (or contraction ratio) of the similitude $f$.
\end{definition}

\begin{definition}[\textcolor{red}{Iterated Function System (IFS)}]
\index{iterated function system (IFS)}
    An \textbf{Iterated Function System (IFS)} is a finite collection of contraction mappings (or similitudes) on a complete metric space $(X, d)$.

\noindent Specifically, an IFS is the collection $\mathcal{I}$ given by:
$$
\mathcal{I} = \{f_1, f_2, \dots, f_N\}
$$
where each function $f_i: X \to X$ is a contraction mapping, meaning there exists a contraction factor $r_i \in [0, 1)$ such that for all $x, y \in X$:
$$
d(f_i(x), f_i(y)) \leq r_i \cdot d(x, y)
$$
The unique non-empty compact set $A \subset X$ that satisfies $A = \bigcup_{i=1}^{N} f_i(A)$ is the \textbf{attractor} of the IFS, which is a self-similar set (see next definition).
\end{definition}

\begin{definition}[\textcolor{red}{Self-Similar Set}]
\index{self-similar set}
A non-empty compact set $A \subset \mathbb{R}^n$ is a \textbf{self-similar set} if it is the unique set that satisfies the equation for a finite collection of similitudes $\mathcal{I} = \{f_1, f_2, \dots, f_N\}$:
$$
A = \bigcup_{i=1}^{N} f_i(A)
$$
\end{definition}

\begin{definition}[\textcolor{red}{Similarity Dimension ($\dim_{\text{sim}}(A)$)}]
\index{similarity dimension}
    
The \textbf{Similarity Dimension} of the set $A$, generated by the IFS $\mathcal{I}$ with scaling ratios $r_1, r_2, \dots, r_N$, is the unique number $s$ that satisfies the characteristic equation:
\begin{equation}
\sum_{i=1}^{N}r_{i}^{s}=1
\end{equation}
If all ratios are equal, $r_i = r$ for all $i$, the formula simplifies to:
\begin{equation}
    \label{eq:simil-dim}
\dim_{\text{sim}}(A)=\frac{\log N}{\log(1/r)}
\end{equation}

\end{definition}
The equality between the similarity and Hausdorff dimension requires a condition known as the Open Set Condition (OSC).
\begin{definition}[\textcolor{red}{Open Set Condition (OSC)}]
\index{open set condition}
An IFS $\mathcal{I} = \{f_1, \dots, f_N\}$ satisfies the \textbf{Open Set Condition (OSC)} if there exists a non-empty, bounded, open set $V \subset \mathbb{R}^n$ such that:
\begin{enumerate}
    \item $f_i(V) \subset V$ for all $i=1, \dots, N$.
    \item The image sets $\{f_i(V)\}_{i=1}^N$ are \textbf{pairwise disjoint}: $f_i(V) \cap f_j(V) = \emptyset$ for all $i \neq j$.
\end{enumerate}    
\end{definition}
The theorem that tells us when the similarity and Hausdorff dimension ar the same is the Moran/Hutchinson theorem.
\begin{theorem}[\textcolor{red}{Moran-Hutchinson Theorem}]
\index{Moran-Hutchinson theorem}
If a self-similar set $A$ is the attractor of an IFS of similitudes $\mathcal{I}=\{f_1, \dots, f_N\}$ and the system satisfies the \textbf{Open Set Condition (OSC)}, then the Hausdorff dimension of $A$ is exactly equal to its Similarity Dimension:

\begin{equation}
\dim_{\mathcal{H}}(A) = \dim_{\operatorname{sim}}(A)
\end{equation}
See for more details \citep{hutchinson_fractals_1981}.
\end{theorem}

\subsection{Square Dimension with the Similarity Method}
\index{similarity method!square}

To apply the method discussed before to $S=[0,1]^2\in \mathbb{R}^2$ we need to find an IFS that satisfy the open set condition. Now naturally a square can be splitted in 4 smaller squares, each with a side length 1/2 of the original one. So we have four $f_1, f_2, f_3$ and $f_4$ each with a scaling ratio $r=1/2$. For example you can imagine $f_1:[0,1|\times[0,1]\rightarrow [0,1/2|\times [0,1/2]$.
It is clearly true that
\begin{equation}
    S=\bigcup_{i=1}^4 f_i(S)
\end{equation}
Now we need to verify the open set condition. We can chose the open unit square $V=(0,1)\times (0,1)$. $f_i$ maps $V$ into a smaller square, and since the set is open, it maps an open set into an open set contained in $V$. Additionally the $f_i(V)$ are clearly disjoint by construction. Thus the OSC condition is verified. We can then use Equation (\ref{eq:simil-dim}) to calculate $\dim_{\text{sim}}(S)$
\begin{equation}
    \dim_{\text{sim}}(S) = \frac{\log 4}{\log 2} = 2
\end{equation}
and, thanks to the Moran-Hutchinson theorem we can say that $\dim_\mathcal{H}(S)=2$.

\subsection{Cantor Set Dimension with the Similarity Method}
\index{similarity method!Cantor set}

As for the square, we need to find an IFS that satisfies the open set condition. At stage $N$, $C_n$ consists of $2^N$ copies of similar sets, so we can define the IFS as $2^N$ $f_i$ with a scaling ratio $r=(1/3)^N$.
Now it is true that
\begin{equation}
    C_N=\bigcup_{i=1}^{2^N} f_i(C_N)
\end{equation}
Now we need to verify the open set condition. We can simply choose the open unit interval $V=(0,1)$. $f_i$ maps $V$ into a smaller subset of $V$ (thus contained in $V$) and additionally by construction the $f_i(V)$ are diskoint. Thus, the OSC condition is verified. We can then use Equation 
(\ref{eq:simil-dim}) to calculate $\dim_{\text{sim}}(C)$ in the limit of $N\rightarrow \infty$.
\begin{equation}
    \dim_{\text{sim}}(C) = \frac{\log 2^N}{\log 3^N} = \frac{\log 2}{\log 3}
\end{equation}
Thus we have that 
\begin{equation}
    \dim_\mathcal{H}(C)=\frac{\log 2}{\log 3}.
\end{equation}


\appendix


\section{Additional Definitions}

An additional definition of $\mathcal{H}_\delta^s$ is given by \cite{evans2015measure}. They define the Hausdorff measure as
\begin{definition}[\textcolor{red}{Hausdorff Measure (Evans and Gariepy)}]
\cite{evans2015measure} define, for example, the Hausdorff measure as the following. Given $A\subset \mathbb{R}^n$, $0\leq s <\infty$, 
    $0<\delta \leq \infty$, they define
    \begin{equation}
        \mathcal{H}_\delta^s(A) = \inf \left\{
            \sum_{j=1}^\infty \alpha(s) 
            \left( \frac{\operatorname{diam } C_j}{2}\right)^s
            | \  A\subseteq\bigcup_{j=1}^\infty C_j, \operatorname{diam } C_j\leq\delta
        \right\}
    \end{equation}
    with
    \begin{equation}
        \alpha(s)\equiv \frac{\pi^{s/2}}{\Gamma(s/2+1)} 
    \end{equation}
        where $\Gamma(s)=\int_0^\infty e^{-x}x^{s-1}$ ($0<s<\infty)$ is the gamma function.
\end{definition}

\begin{definition}[\textcolor{red}{Covering (from \cite{lio_analysis_nodate}}]
\label{eq:covering}
    This definition is given in \citep{lio_analysis_nodate}. Let $\mathcal{K}\subseteq\mathcal{P}(X)$. $\mathcal{K}$ is called a \textbf{covering}\index{covering} of $X$ if
    \begin{enumerate}
        \item $\emptyset \in \mathcal{K}$
        \item $\exists (K_j)_{j\in \mathbb{N}}\subseteq \mathcal{K}: X=\bigcup_{j=1}^\infty K_j$
    \end{enumerate}
\end{definition}

\begin{definition}[\textcolor{red}{Covering (less restrictive version)}]
\label{def:covering2}
    A covering of a set $E$ is a collection of sets $(C_i)_{i\in I}$ ($I$ at most countable set of indeces) such that
    \begin{equation}
        E\subseteq \bigcup_{i=1}^\infty C_i
    \end{equation}
\end{definition}

\bibliographystyle{plainnat} 
\bibliography{references}

@book{evans2015measure,
  title={Measure theory and fine properties of functions, revised edition},
  author={Evans, Lawrence C. and Gariepy, Ronald F.},
  year={2015},
  publisher={Chapman and Hall/CRC},
  isbn={9781003583004}
}

@book{lio_analysis_nodate,
	title = {Analysis 3 ({Measure} {Theory}) {Fall} {Semester} 2025 Lecture Notes},
	language = {en},
	author = {Da Lio, Francesca},
    publisher = {ETH, Z\"urich},
year ={2025}
}

@book{falconer2003fractal,
  title     = {Fractal Geometry: Mathematical Foundations and Applications},
  author    = {Falconer, Kenneth},
  year      = {2003},
  publisher = {John Wiley \& Sons},
  address   = {Chichester},
  isbn      = {9780470848616},
  doi       = {10.1002/0470013850}
}

@article{hutchinson_fractals_1981,
	title = {Fractals and {Self} {Similarity}},
	volume = {30},
	issn = {0022-2518},
	url = {https://www.jstor.org/stable/24893080},
	number = {5},
	urldate = {2025-11-12},
	journal = {Indiana University Mathematics Journal},
	author = {Hutchinson, John E.},
	year = {1981},
	note = {Publisher: Indiana University Mathematics Department},
	pages = {713--747}
}

@article{hao_hausdorff_nodate,
    author ={Wanqing Hao},
	title = {{The} {Hausdorff} {dimension}: {construction} {and} {methods} {calculation} {calculation}},
    year = {2024},
    url={https://math.uchicago.edu/~may/REU2024/REUPapers/Hao,Wanqing(Anita).pdf}
}

@article{cao_several_2023,
	title = {Several {Methods} for {Calculating} and {Estimating} the {Hausdorff} {Dimension}},
	volume = {5},
	copyright = {https://creativecommons.org/licenses/by/4.0},
	issn = {2832-6024},
	url = {https://drpress.org/ojs/index.php/fcis/article/view/11588},
	doi = {10.54097/fcis.v5i1.11588},
	language = {en},
	number = {1},
	urldate = {2025-11-12},
	journal = {Frontiers in Computing and Intelligent Systems},
	author = {Cao, Ziang},
	month = sep,
	year = {2023},
	pages = {28--32},
}

\printindex
\end{document}